\documentclass{gtart}

\def\ifplaintex{\expandafter\ifx\csname documentclass\endcsname\relax}

\def\gtp{{\mathsurround=0pt\it $\cal G\mskip-2mu$eometry \&\ 
$\cal T\!\!$opology $\cal P\!$ublications}}  

\def\recd{{\small Received:\qua\receiveddate\ifx\reviseddate\relax
\else\qquad Revised:\qua\reviseddate\fi\par}} 

\def\lognumber#1{\def\thelognumber{#1}}
\def\volumenumber#1{\def\thevolumenumber{#1}}
\def\volumeyear#1{\def\thevolumeyear{#1}}
\def\papernumber#1{\def\thepapernumber{#1}}
\def\pagenumbers#1#2{\def\startpage{#1}\def\finishpage{#2}}
\def\published#1{\def\publishdate{#1}}

\def\received#1{\def\receiveddate{#1}}
\def\revised#1{\def\reviseddate{#1}}
\def\accepted#1{\def\accepteddate{#1}}

\long\def\asciiabstract#1{\long\def\theasciiabstract{#1}}

\let\\\par\let\thelognumber\relax\let\thevolumenumber\relax
\let\thepapernumber\relax\let\thevolumeyear\relax\let\startpage\relax
\let\finishpage\relax\let\publishdate\relax\let\receiveddate\relax
\let\reviseddate\relax\let\accepteddate\relax\let\theasciititle\relax
\let\theasciiauthors\relax
\let\theasciiabstract\relax

\let\theasciiemail\relax

\ifplaintex
\font\logobig=cmssbx10 scaled 3836
\font\logomed=cmssbx10 scaled 2557
\else
\font\logobig=cmssbx10 scaled 4200
\font\logomed=cmssbx10 scaled 2800
\fi

\long\def\makeagttitle{   
\count0=\startpage
\agt\hfill      
\hbox to 45truept{\vbox to 0pt{\vglue -13truept{\logomed A\kern -.37em{\logobig 
T}\kern -.38em G}\vss}\hss}
\break
{\small Volume \thevolumenumber\ (\thevolumeyear)
\startpage--\finishpage\nl
Published: \publishdate}

\vglue .25truein

{\parskip=0pt\leftskip 0pt plus
1fil\def\\{\par\smallskip}{\Large\bf\thetitle}\par\medskip} \vglue
0.05truein

{\parskip=0pt\leftskip 0pt plus 1fil\def\\{\par}{\sc\theauthors}
\par\medskip}%
 
\vglue 0.03truein 

{\small\leftskip 25truept\rightskip 25truept{\bf Abstract}\stdspace\theabstract

{\bf AMS Classification}\stdspace\theprimaryclass
\ifx\thesecondaryclass\relax\else; \thesecondaryclass\fi\par
{\bf Keywords}\stdspace \thekeywords\par}\vglue 7truept

}   

\ifplaintex
\font\phead=cmsl9 scaled 950
\font\pnum=cmbx10 scaled 913
\font\pfoot=cmsl9 scaled 950
\headline{\vbox to 0pt{\vskip -4.5mm\line{\small\phead\ifnum
\count0=\startpage ISSN 1472-2739 (on-line) 1472-2747 (printed)
\hfill {\pnum\folio}\else\ifodd\count0\def\\{ }%
\ifx\theshorttitle\relax\thetitle\else\theshorttitle\fi\hfill{\pnum\folio}
\else\def\\{ and }{\pnum\folio}\hfill\ifx\theshortauthors\relax\theauthors
\else\theshortauthors\fi\fi\fi}\vss}}
\footline{\vbox to 0pt{\vglue 0mm\line{\small\pfoot\ifnum\count0=\startpage
\copyright\ \gtp\hfill\else
\agt, Volume \thevolumenumber\ (\thevolumeyear)\hfill\fi}\vss}}
\else
\font=cmsl9 scaled 1050
\font\lnum=cmbx10 
\font=cmsl9 scaled 1050
\makeatletter
\def\@oddhead{{\small\lhead\ifnum\count0=\startpage ISSN 1472-2739 
(on-line) 1472-2747 (printed)\hfill {\lnum\number\count0}\else\ifodd\count0
\def\\{ }\ifx\theshorttitle\relax \thetitle \else\theshorttitle\fi\hfill
{\lnum\number\count0}\else\def\\{ and }{\lnum\number\count0}
\hfill\ifx\theshortauthors\relax 
\theauthors\else\theshortauthors\fi\fi\fi}}\def\@evenhead{\@oddhead}
\def\@oddfoot{\small\lfoot\ifnum\count0=\startpage\copyright\ \gtp\hfill\else
\agt, Volume \thevolumenumber\ (\thevolumeyear)\hfill\fi}
\def\@evenfoot{\@oddfoot}
\makeatother
\fi
\let\maketitlepage\makeagttitle
\let\makeshorttitle\maketitlepage
\let\maketitle\maketitlepage

\newwrite\gtoutfile
\long\gdef\makeheadfile{  
{\def\\{, }\def\s{ }
\immediate\openout\gtoutfile head.xxx
\immediate\write\gtoutfile{To: math@arxiv.org}
\immediate\write\gtoutfile{Subject: put OR rep NNNNN:ppppp}
\immediate\write\gtoutfile{--text follows this line--}
\immediate\write\gtoutfile{Proxy-for: \ifx\theasciiauthors\relax
\theauthors\else\theasciiauthors\fi\s<\ifx\theasciiemail\relax\theemail\else\theasciiemail\fi>}
\immediate\write\gtoutfile{\noexpand\\}
\immediate\write\gtoutfile{Authors: \ifx\theasciiauthors\relax
\theauthors\else\theasciiauthors\fi}
{\def\\{ }\immediate\write\gtoutfile{Title: \ifx\theasciititle\relax
\thetitle\else\theasciititle\fi}}
\immediate\write\gtoutfile{Subj-class: GT or SG, GR etc}
\immediate\write\gtoutfile{MSC-class: \theprimaryclass\ifx\thesecondaryclass\relax\else, \thesecondaryclass\fi}
\immediate\write\gtoutfile{Journal-ref: Algebr. Geom. Topol. \thevolumenumber\s
(\thevolumeyear) \startpage-\finishpage}
\immediate\write\gtoutfile{Comments: Published by Algebraic and
Geometric Topology at}
\immediate\write\gtoutfile{\s\s\s  http://www.maths.warwick.ac.uk/agt/AGTVol\thevolumenumber/agt-\thevolumenumber-\thepapernumber.abs.html}
\immediate\write\gtoutfile{\noexpand\\}
\immediate\write\gtoutfile{}
\ifx\theasciiabstract\relax
\immediate\write\gtoutfile{\theabstract}\else
\immediate\write\gtoutfile{\theasciiabstract}\fi
\immediate\write\gtoutfile{}
\immediate\write\gtoutfile{\noexpand\\}
\immediate\write\gtoutfile{}
\immediate\closeout\gtoutfile}}  

\def\maketitlepage{\makeagttitle\makeheadfile}
\let\makeshorttitle\maketitlepage
\let\maketitle\maketitlepage

\lognumber{20}
\volumenumber{3}
\volumeyear{2003}
\papernumber{20}
\published{21 June 2003}
\pagenumbers{587}{591}
\received{18 August 2002}
\revised{15 June 2003}
\accepted{23 October 2002}

\usepackage{amsmath,pstricks,subfigure}

\newtheorem{theorem}{Theorem}

\theoremstyle{remark}

\newcommand{\sm}{\setminus}

\newcommand{\ie}{\textit{i.e.}}

\newenvironment{fullfigure}[2]
    {\begin{figure}[htb]\begin{center}\def\ffa{#1}\def\ffb{#2}}
    {\caption{\ffb.}\label{\ffa}\end{center}\end{figure}}

\newcommand{\fig}[1]{Figure~\ref{#1}}
\newcommand{\thm}[1]{Theorem~\ref{#1}}

\psset{linewidth=.5pt,dash=2.5pt 2.5pt,unit=.2in,arrowsize=3pt 3}
\SpecialCoor

\begin{document}
\title{What is a virtual link?}
\author{Greg Kuperberg}
\address{Department of Mathematics,
    University of California\\Davis, CA 95616, USA}
\email{greg@math.ucdavis.edu}

\begin{abstract}
Several authors have recently studied virtual knots and links because they
admit invariants arising from $R$-matrices.  We prove that every virtual link
is uniquely represented by a link $L \subset S \times I$ in a thickened,
compact, oriented surface $S$ such that the link complement $(S \times I) \sm
L$ has no essential vertical cylinder.
\end{abstract}
\begin{abstract}
Several authors have recently studied virtual knots and links because they
admit invariants arising from $R$-matrices.  We prove that every virtual link
is uniquely represented by a link $L \subset S \times I$ in a thickened,
compact, oriented surface $S$ such that the link complement $(S \times I) \sm
L$ has no essential vertical cylinder.
\end{abstract}
\asciiabstract{Several authors have recently studied virtual knots
and links because they admit invariants arising from R-matrices.  We
prove that every virtual link is uniquely represented by a link L
in S X I, a thickened, compact, oriented surface S, such
that the link complement (S X I) - L  has no essential
vertical cylinder.}


\maketitle

A \emph{virtual link} $L$ is an equivalence class of decorated, finite,
tetravalent graphs $\Gamma$.  The edges at each vertex must be cyclically
ordered, and two opposite edges are labelled as an \emph{overcrossing}, while
the other two are labelled as an \emph{undercrossing}.  The equivalence
relation is the one given by Reidemeister moves. The notion was proposed by
Kauffman \cite{Kauffman:virtual} in light of the fact that $R$-matrices and
quandles, which are commonly used to make link invariants, also yield
invariants of virtual links.

We will borrow from the analysis of virtual links by Carter, Kamada, and Saito
\cite{CKS:stable}.  They show that virtual links are equivalent to stable
equivalence classes of links projections onto compact, oriented surfaces.  
(Fenn, Rourke, and Sanderson have written a self-contained proof
\cite{FRS:rack}.) The surface $S$ need not be connected, but we require that
each component contains at least one component of the projection $P$. The
projection $P$ is considered up to Reidemeister moves, and the stabilization
operation consists of adding a handle to $S$. Note that the feet of a
stabilizing handle may lie in any two regions of the complement of $P$
(\fig{f:stab}).

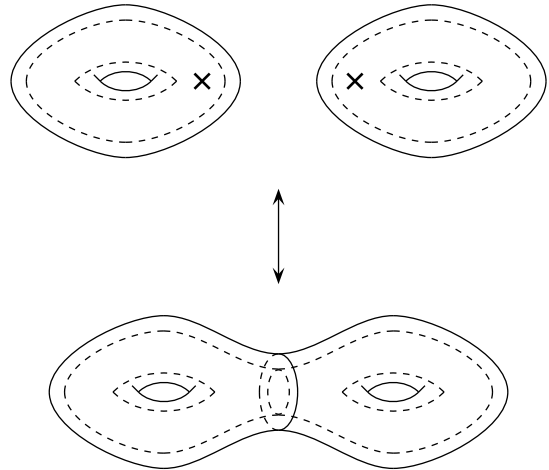
\begin{fullfigure}{f:stab}{The stabilization move on $S \times I$}
\pspicture(-7.5,-2.5)(7.5,2.5)
\psbezier(1,0)(1,1)(3,2)(4,2) \psbezier(4,2)(5,2)(7,1)(7,0)
\psbezier(1,0)(1,-1)(3,-2)(4,-2) \psbezier(4,-2)(5,-2)(7,-1)(7,0)
\psbezier[linestyle=dashed](1.4,0)(1.4,.8)(3,1.6)(4,1.6)
\psbezier[linestyle=dashed](4,1.6)(5,1.6)(6.6,.8)(6.6,0)
\psbezier[linestyle=dashed](1.4,0)(1.4,-.8)(3,-1.6)(4,-1.6)
\psbezier[linestyle=dashed](4,-1.6)(5,-1.6)(6.6,-.8)(6.6,0)
\psbezier(3.333,0)(3.667,.333)(4.333,.333)(4.667,0)
\psbezier(4.667,0)(4.333,-.333)(3.667,-.333)(3.333,0)
\psline(3.333,0)(3.167,.167) \psline(4.667,0)(4.833,.167)
\psbezier[linestyle=dashed](2.667,0)(3.333,.667)(4.667,.667)(5.333,0)
\psbezier[linestyle=dashed](5.333,0)(4.667,-.667)(3.333,-.667)(2.667,0)
\psline[linewidth=1pt](1.8,.2)(2.2,-.2)\psline[linewidth=1pt](1.8,-.2)(2.2,.2)
\psbezier(-1,0)(-1,1)(-3,2)(-4,2) \psbezier(-4,2)(-5,2)(-7,1)(-7,0)
\psbezier(-1,0)(-1,-1)(-3,-2)(-4,-2) \psbezier(-4,-2)(-5,-2)(-7,-1)(-7,0)
\psbezier[linestyle=dashed](-1.4,0)(-1.4,.8)(-3,1.6)(-4,1.6)
\psbezier[linestyle=dashed](-4,1.6)(-5,1.6)(-6.6,.8)(-6.6,0)
\psbezier[linestyle=dashed](-1.4,0)(-1.4,-.8)(-3,-1.6)(-4,-1.6)
\psbezier[linestyle=dashed](-4,-1.6)(-5,-1.6)(-6.6,-.8)(-6.6,0)
\psbezier(-3.333,0)(-3.667,.333)(-4.333,.333)(-4.667,0)
\psbezier(-4.667,0)(-4.333,-.333)(-3.667,-.333)(-3.333,0)
\psline(-3.333,0)(-3.167,.167) \psline(-4.667,0)(-4.833,.167)
\psbezier[linestyle=dashed](-2.667,0)(-3.333,.667)(-4.667,.667)(-5.333,0)
\psbezier[linestyle=dashed](-5.333,0)(-4.667,-.667)(-3.333,-.667)(-2.667,0)
\psline[linewidth=1pt](-1.8,.2)(-2.2,-.2)
\psline[linewidth=1pt](-1.8,-.2)(-2.2,.2)
\endpspicture \\
\pspicture(-1,-1.5)(1,1.5) \psline{<->}(0,1.25)(0,-1.25) \endpspicture \\
\pspicture(-6.5,-2.5)(6.5,2.5)
\psellipse[linestyle=dashed](0,0)(.3,.6)
\psbezier(.5,0)(.5,.5523)(.2761,1)(0,1)
\psbezier(.5,0)(.5,-.5523)(.2761,-1)(0,-1)
\psbezier[linestyle=dashed](-.5,0)(-.5,.5523)(-.2761,1)(0,1)
\psbezier[linestyle=dashed](-.5,0)(-.5,-.5523)(-.2761,-1)(0,-1)
\psbezier(0,1)(1,1)(2,2)(3,2) \psbezier(3,2)(4,2)(6,1)(6,0)
\psbezier(6,0)(6,-1)(4,-2)(3,-2) \psbezier(3,-2)(2,-2)(1,-1)(0,-1)
\psbezier[linestyle=dashed](0,.6)(1,.6)(2,1.6)(3,1.6)
\psbezier[linestyle=dashed](3,1.6)(4,1.6)(5.6,.8)(5.6,0)
\psbezier[linestyle=dashed](5.6,0)(5.6,-.8)(4,-1.6)(3,-1.6)
\psbezier[linestyle=dashed](3,-1.6)(2,-1.6)(1,-.6)(0,-.6)
\psbezier(2.333,0)(2.667,.333)(3.333,.333)(3.667,0)
\psbezier(3.667,0)(3.333,-.333)(2.667,-.333)(2.333,0)
\psline(2.333,0)(2.167,.167) \psline(3.667,0)(3.833,.167)
\psbezier[linestyle=dashed](1.667,0)(2.333,.667)(3.667,.667)(4.333,0)
\psbezier[linestyle=dashed](4.333,0)(3.667,-.667)(2.333,-.667)(1.667,0)
\psbezier(0,1)(-1,1)(-2,2)(-3,2) \psbezier(-3,2)(-4,2)(-6,1)(-6,0)
\psbezier(-6,0)(-6,-1)(-4,-2)(-3,-2) \psbezier(-3,-2)(-2,-2)(-1,-1)(0,-1)
\psbezier[linestyle=dashed](0,.6)(-1,.6)(-2,1.6)(-3,1.6)
\psbezier[linestyle=dashed](-3,1.6)(-4,1.6)(-5.6,.8)(-5.6,0)
\psbezier[linestyle=dashed](-5.6,0)(-5.6,-.8)(-4,-1.6)(-3,-1.6)
\psbezier[linestyle=dashed](-3,-1.6)(-2,-1.6)(-1,-.6)(0,-.6)
\psbezier(-2.333,0)(-2.667,.333)(-3.333,.333)(-3.667,0)
\psbezier(-3.667,0)(-3.333,-.333)(-2.667,-.333)(-2.333,0)
\psline(-2.333,0)(-2.167,.167)
\psline(-3.667,0)(-3.833,.167)
\psbezier[linestyle=dashed](-1.667,0)(-2.333,.667)(-3.667,.667)(-4.333,0)
\psbezier[linestyle=dashed](-4.333,0)(-3.667,-.667)(-2.333,-.667)(-1.667,0)
\endpspicture
\end{fullfigure}

As \fig{f:stab} also shows, the reverse destabilization operation consists of
cutting the surface $S$ along a circle $C$ which is disjoint from $P$, and
capping the resulting boundary. If $C$ separates the connected component of
$S'$ of $S$ containing it, then we require that both components of $S' \sm C$
contain part of $P$; otherwise destabilization would create a naked surface
component.

It is also well-known that a link drawn on a surface $S$, considered up to
Reidemeister moves, is equivalent to a (tame) link $L \subset S \times I$ in
the thickened surface $S$, considered up to isotopy.  The destabilization
operation, then, consists of cutting $S$ along a vertical annulus $C \times I$
which is disjoint from $L$, and capping the two resulting annuli with thickened
disks $D^2 \times I$.  Since $L$ is only considered up to isotopy, it is
equivalent to allow this operation with any topologically vertical annulus $A$,
\ie, any properly embedded annulus isotopic to some $C \times I$.  This will be
our working definition of virtual links.

\begin{theorem} Every stable equivalence class of links in thickened surfaces
has a unique irreducible representative.
\label{th:main}
\end{theorem}

\thm{th:main} is proved in the same way as several classical results in
3-manifold topology, from the unique factorization of knots and 3-manifolds
\cite{Schubert:knotens,Milnor:unique} to the Jaco-Shalen-Johannson theorem
\cite{JS:seifert,Johannson:homotopy}, from which \thm{th:main} also follows as
a corollary.  Another similar result is the author's classification of scalar
involutory Hopf words \cite{Kuperberg:hopf}. The motivation is also
similar, since a virtual link can be viewed as a scalar $R$-matrix word.

\begin{proof}  In outline, we induct on the complexity of intersection between
two destabilization annuli by compressing one along an innermost disk of the
other.

First generalize the definition of destabilization of a link $L \subset S
\times I$ to include two other operations: 

\textbf{(1)}\qua If a sublink $L'
\subset L$ is separated from the rest of $L$ by a sphere or disk $A \subset S
\times I$, then we can remove $L'$ from $S \times I$ and place it in a separate
thickened sphere $S^2 \times I$. 

\textbf{(2)}\qua If an annulus $A$ divides $S
\times I$ with $L$ entirely on one side and some genus of $S$ on the other
side, we can cut $S \times I$ along $A$, discard the naked component, and cap
the remaining component. Both operations can easily be reproduced by
destabilization as defined previously.

Say that a surface is \emph{admissible} if it is a vertical annulus, a sphere,
or a proper disk; and that an  admissible surface is \emph{essential} if it
does not bound a ball in $(S \times I) \sm L$.  Thus, admissible, essential
surfaces are those along which we can destabilize $L \subset S \times I$.

Suppose, to the contrary of the conclusion, that some link $L \subset S \times
I$ has more than one irreducible descendant.  If $S$ has $c$ components with
total genus $g$, and $L$ has $n$ components, assume that $g+n-c$ is minimal
among counterexamples.  (Note that $n\ge c$.)  Then every destabilization of $L
\subset S \times I$ has a unique irreducible descendant, since destabilization
always reduces $g+n-c$. Say that two such destabilizations, $L \subset S_1
\times I$ and $L \subset S_2 \times I$, are \emph{descent equivalent} if their
irreducible descendants are isomorphic. The aim is to show that all
destabilizations of $L \subset S \times I$ are descent equivalent.

For example, if $A_1$ and $A_2$ are disjoint admissible, essential surfaces,
then the resulting destabilizations $L \subset S_1 \times I$ and $L \subset S_2
\times I$ are descent equivalent.  This is immediate if $A_1$ and $A_2$ are
parallel.  If they are not, then we can destabilize each $L \subset S_i \times
I$ along $A_{3-i}$ to produce a common descendant.

Suppose that $A_1$ and $A_2$ are descent-inequivalent surfaces in general
position, and that they intersect in the fewest curves among
descent-inequivalent pairs in $(S \times I) \sm L$.   If a curve $C \subset A_1
\cap A_2$ is a circle, then it is either \emph{horizontal} in $A_i$, if $A_i$
is an annulus and $C$ is parallel to $\partial A_i$, or it bounds a disk in
$A_i$.  If $C$ is an arc, then it is either \emph{vertical}, if $A_i$ is an
annulus and $C$ connects the two components of $\partial A_i$, or it and part
of the boundary of $A_i$ bound a disk.

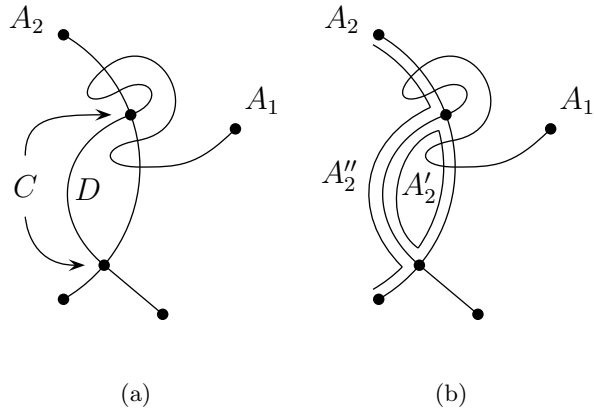
\begin{fullfigure}{f:innermost}{Compressing $A_2$ along the disk $D$ to
    simplify $A_1 \cap A_2$ (side view)}
\subfigure[]{\pspicture(0,-4.5)(8,5)
\qdisk(4;-60){.15} \qdisk(4;60){.15} \qdisk(4;-40){.15} \qdisk(6;-40){.15}
\psarc(0,0){4}{-60}{60} \psline(4;-40)(6;-40)
\psbezier(4;-40)(2;-40)(2;20)(4;20) \qdisk(4;20){.15}
\psbezier(4;20)(5;20)(5;30)(4;30) \psbezier(4;30)(3;30)(3;40)(4;40)
\psbezier(4;40)(6;40)(6;10)(4;10) \psbezier(4;10)(3;10)(3;0)(4;0)
\psbezier(4;0)(5;0)(5.5,0)(6.5,1) \qdisk(6.5,1){.15}
\rput(4;75){\large $A_2$} \rput(7.2,1.75){\large $A_1$}
\rput(2.7;-12.5){\large $D$} \rput(1,-.5){\large $C$}
\psbezier{->}(1,.3)(1,1.3)(1.759,1.368)(3.259,1.368)
\psbezier{->}(1,-1.3)(1,-1.3)(1.064,-2.571)(2.564,-2.571)
\endpspicture}
\subfigure[]{\pspicture(0,-4.5)(8,5)
\qdisk(4;-60){.15} \qdisk(4;60){.15} \qdisk(4;-40){.15} \qdisk(6;-40){.15}
\psarc(0,0){4}{-60}{60} \psline(4;-40)(6;-40)
\psbezier(4;-40)(2;-40)(2;20)(4;20) \qdisk(4;20){.15}
\psbezier(4;20)(5;20)(5;30)(4;30) \psbezier(4;30)(3;30)(3;40)(4;40)
\psbezier(4;40)(6;40)(6;10)(4;10) \psbezier(4;10)(3;10)(3;0)(4;0)
\psbezier(4;0)(5;0)(5.5,0)(6.5,1) \qdisk(6.5,1){.15}
\psarc(0,0){3.7}{-35}{15}
\psbezier(3.7;15)(2.5;15)(2.5;-35)(3.7;-35)
\psarc(0,0){3.7}{-60}{-45}
\psbezier(3.7;-45)(1.7;-45)(1.7;25)(3.7;25)
\rput(4;75){\large $A_2$} \rput(7.2,1.75){\large $A_1$}
\rput(3.1;-10){\large $A'_2$} \rput(1;-13){\large $A''_2$}
\psarc(0,0){3.7}{25}{60}
\endpspicture}
\end{fullfigure}

If a circle of $A_1 \cap A_2$ is non-horizontal in $A_1$ (say), then some such
circle $C$ is \emph{innermost}, meaning that it bounds a naked disk $D$ in
$A_1$, as in \fig{f:innermost}(a). In this case let $A'_2$ and $A''_2$ be the
connected components of the compression of $A_2$ along the disk $D$, as in
\fig{f:innermost}(b).  Both $A'_2$ and $A''_2$ are admissible, and at least one
is essential, for otherwise $A_2$ would not be.  If $A'_2$ (say) is essential,
then it intersects $A_1$ less than $A_2$ does, and it does not intersect $A_2$
at all.  But since $A_1$ and $A_2$ intersect least among descent-inequivalent
pairs of essential surfaces, it would follow that $A_1$ and $A_2$ are both
descent-equivalent to $A'_2$, a contradiction.

The same argument applies if $A_1 \cap A_2$ has a non-vertical arc in $A_1$,
and of course also to non-horizontal circles and non-vertical arcs in $A_2$.
Thus $A_1 \cap A_2$ consists entirely of vertical segments or horizontal
circles in both $A_1$ and $A_2$.  In particular, both $A_1$ and $A_2$ are
vertical annuli and not disks or spheres.

Suppose that $A_1 \cap A_2$ consists of horizontal circles. We can assume that
none of the four circles of $\partial A_1$ and $\partial A_2$ bounds a disk in
$S \times \partial I$. If, say, $C \subset A_1$ bounds a naked disk $D \subset
S \times \partial I$, then we replace $A_1$ by the disk $D \cup A_1$ and reduce
to a previous case without worsening $A_1 \cap A_2$. Otherwise let $C \subset
A_1 \cap A_2$ be an \emph{outermost} circle, meaning that it and one component
of $\partial A_1$ bound a naked annulus $A \subset A_1$.  The circle $C$
divides $A_2$ into two annuli $A'_2$ and $A''_2$, one of which, say $A'_2$,
makes a vertical annulus together with $A$. The annulus $A \cup A'_2$ is
necessarily essential since its boundary circles do not bound disks in $S
\times \partial I$.  But after displacement, $A \cup A'_2$ intersects $A_1$ and
$A_2$ less than they do each other.  It is therefore descent equivalent to
both.

Finally suppose that $A_1 \cap A_2$ consists of vertical arcs.  The boundary of
a regular neighborhood of $A_1 \cup A_2$ consists of vertical annuli
$B_1,B_2,\ldots,B_n$.  Each $B_i$ is disjoint from both $A_1$ and $A_2$, so if
any of them is essential, it is descent equivalent to both $A_1$ and $A_2$, a
contradiction.  But if they are all inessential, then one of them, say $B_1$,
separates $A_1$ and $A_2$ from the link $L$ and bounds a ball that contains
$A_1$ and $A_2$. This contradicts the hypothesis that $A_1$ and $A_2$ are
essential.
\end{proof}

\thm{th:main} implies that if two links are equivalent as virtual links, then
they are equivalent as links. It also generalizes to oriented virtual links, to
colored virtual links, and even to virtual tangled graphs.  The proof in each
case is the same.

\subsection*{Acknowledgments}

The author would like to thank Dror Bar-Natan, Louis Kauffman, Colin Rourke,
and Dylan Thurston for bringing the virtual link problem to his attention, and
for very helpful discussions.

The author was supported by NSF grant DMS \#0072342.

\Addresses\recd

\end{document}